# Universal Malliavin Calculus in Fock and Lévy-Itô Spaces

David Applebaum,
Probability and Statistics Department,
University of Sheffield,
Hicks Building, Hounsfield Road,
Sheffield, England, S3 7RH

e-mail: D.Applebaum@sheffield.ac.uk

**Abstract**

We review and extend Lindsay's work on abstract gradient and divergence operators in Fock space over a general complex Hilbert space. Precise expressions for the domains are given, the $L^2$-equivalence of norms is proved and an abstract version of the Itô-Skorohod isometry is established. We then outline a new proof of Itô's chaos expansion of complex Lévy-Itô space in terms of multiple Wiener-Lévy integrals based on Brownian motion and a compensated Poisson random measure. The duality transform now identifies Lévy-Itô space as a Fock space. We can then easily obtain key properties of the gradient and divergence of a general Lévy process. In particular we establish maximal domains of these operators and obtain the Itô-Skorohod isometry on its maximal domain.

# 1 Introduction

Malliavin calculus is one of the deepest and most important areas within contemporary stochastic analysis. It was originally developed as a new probabilistic technique to find smooth densities for solutions of stochastic differential equations (SDEs). At a more fundamental level it provides an intrinsic differential calculus in Gaussian probability spaces based on two mutually adjoint linear operators - the gradient and the divergence (see e.g. [36], [39], [27], [23], [50] for monograph accounts). More recently it has enabled the developments of new techniques in mathematical finance (see [37] and references therein).

Ever since the early days of the subject there has been plenty of activity in widening the scope of Malliavin calculus to include jump processes and [8] is a monograph dedicated to this theme. More recently there has been increased interest in these ideas - partly due to new progress in finding smooth densities for special classes of SDEs driven by Lévy processes (see e.g. [28], [24]) but also for the need to extend the calculus to financial models based on jump processes (see e.g. [31], [11], [44] and the forthcoming monograph [15]).

Fock space has long been known to be intimately connected with probability theory. Indeed if the natural $L^2$-space of a process has a chaotic decomposition then it is automatically isomorphic to a Fock space over the Cameron-Martin space (in probabilistic language) or one-particle space (in physical terms). The first direct use of Fock space ideas in Malliavin calculus seems to have been by Dermoune, Krée and Wu [13] in work on non-anticipating stochastic calculus (including a generalised Itô formula) for the Poisson process. These ideas were then taken up by Nualart and Vives who defined the gradient and divergence for the Poisson process directly in Fock space [40] and Dermoune [12] who extended the work of [13] to general Lévy processes.

In a separate development, Hudson and Parthasarathy [22] realised that Fock space is the natural setting for a quantum stochastic calculus based on a non-commutative splitting of Brownian motion and the Poisson process into constituent annihilation, creation and conservation noise processes (see also [43], [38], [34]). Parthasarathy (see [43] p.155-8) also showed that Lévy processes may be represented in Fock space and the corresponding extended quantum stochastic calculus was developed in [1], [3]. A key paper by Lindsay [33] brought Malliavin calculus directly into the non-commutative framework to devise a non-anticipating quantum stochastic calculus. Similar ideas appear in the works of Belavkin (see e.g. [5]). More recent developments in this area can be found in [4]. A fully quantised Malliavin calculus based on the Wigner density is due to Franz, Léandre and Schott [19], [20].

The goal in this paper is to begin to develop a *universal Malliavin calculus* in Fock space over a general Hilbert space (see also Privault and Wu [48]). There is no probability content (either classical or quantum) in the theory at this stage. We focus our studies on two operators originally introduced by Lindsay [33] in the context of a Fock-Guichardet space and called abstract gradient and divergence therein. We prefer to call them *universal annihilation and creation operators* as they can be transformed into the usual annihilation and creation operators indexed by a given vector in one-particle space after composition with a suitable Dirac bra or ket operator. Following Privault [47] (see also Privault and Wu [48]) we denote these by $\nabla^-$ and $\nabla^+$ respectively. The aim of universal Malliavin calculus can be summed up succinctly as follows - *given a process having a chaos decomposition, map $\nabla^-$ and $\nabla^+$ unitarily into the $L^2$-space of the process where they become the gradient $D$ and divergence $\delta$. Then structural properties of $\nabla^-$ and $\nabla^+$ are automatically transferred to $D$ and $\delta$ with little additional effort.* In this paper we illustrate this technique through application to a Lévy process, but it could easily be applied to any other process having a chaotic representation, such as the Azéma martingale [18] or the Dunkl process [21].

The results obtained for $\nabla^-$ and $\nabla^+$ in [33] were extended to the general case by Lindsay in Proposition 3.1 of [34]. He established three key properties of these operators:

(1) The maximal domain of $\nabla^-$.

(2) The factorisation of the number operator $N = \nabla^+ \nabla^-$.

(3) An isometry-type property for $\nabla^+$ which generalises the key Itô-Skorohod isometry that lies at the heart of non-anticipating stochastic calculus (a result of this kind was also established independently by Privault and Wu [48]).

The proofs of these results were outlined in [33]. We give full proofs in section 3 of this paper for the sake of completeness. The approach presented here is different and we reformulate the result of (1) in a way that will be more familiar to probabilists. We also extend the theory by obtaining a result on the $L^2$-equivalence of norms which is the first step towards a theory of infinite dimensional Sobolev spaces at this level.

In section 4 we turn to probability theory and study the chaotic representation of a Lévy process. This result, originally due to Itô [25], shows that the natural $L^2$-space $\mathcal{H}$ of the process (called *Lévy-Itô space* herein) is naturally isomorphic to the infinite direct sum of the chaoses generated by multiple Wiener-Lévy integrals constructed from the Brownian motion $B$ and compensated Poisson measure $\tilde{N}$ associated to the process through its Lévy-Itô decomposition. More straightforward proofs of this result have recently been found by Løkka [35] for square integrable processes and Petrou [44] in the general case. Their approach is to iterate the Itô representation of the process whereby any element of $\mathcal{H}$ is a constant plus an Itô stochastic integral with respect to $B$ and $\tilde{N}$. The Itô representation is itself proved by a density argument using a class of exponential martingales of Wiener-Lévy stochastic integrals of deterministic functions. We briefly outline a generalisation of these results. The main difference for us is that $\mathcal{H}$ is complex and this allows a simpler proof of the Itô representation using a more natural class of exponential martingales. We only give outline proofs here as the methodology is well-known and a full account will appear shortly in [2] (ii) (similar ideas are employed in Bichteler [7], p.259.)

The chaotic representation induces a unitary isomorphism called the *duality transform* between $\mathcal{H}$ and a certain Fock space. In section 4 we also prove that the image of the exponential vectors under this isomorphism is the stochastic or Doléans-Dade exponentials. A general result of this type was hinted at by Meyer in [38] p.71 but he only explicitly considered the Wiener space case. These vectors and their chaos expansions have recently found interesting applications to interest-rate modelling [9].

In the last part of the paper we apply the duality transform to the results (1) to (5) of section 1 to obtain the Lévy process versions of these for $D$ and $\delta$. Our main results are

1. The maximal domains of $D$ and $\delta$ expressed in terms of chaos expansions.

2. The full Itô-Skorohod isometry on a maximal domain.

We note that the Itô-Skorohod isometry has also recently been investigated in [14] using white noise analysis techniques but these authors were restricted to using a pure jump square integrable Lévy process without drift and no explicit domain was given. A white noise approach is also developed in [30] but under the constraint that the Lévy measure has moments to all orders (see also [29]).

In this paper we have only made the first few steps in the direction of a universal Malliavin calculus. A defect of the theory as it stands is that it only works at the Hilbert space level and so there are, for example, no direct analogues of the full range of infinite dimensional Sobolev spaces which require $L^p$ structure when $p \neq 2$. It may be that this can be remedied by using a Banach-Fock space as in section 6 of [6].

Some related work on Malliavin calculus for Lévy processes has recently been presented in [51]. The key novel ingredient here is the development of a new canonical construction for pure jump Lévy processes which facilitates the study of the Malliavin derivative through its representation as a difference quotient.

Notation: All inner products in complex Hilbert spaces are conjugate linear on the left. The algebraic tensor product of two vector spaces $V_1$ and $V_2$ will be denoted $V_1 \underline{\otimes} V_2$. If $T$ is a closable operator in a Hilbert space, we will throughout this paper use the same notation $T$ for its closure on the larger domain. $\text{Dom}(T)$ will always denote the maximal domain of $T$. If $S$ is a topological space then $\mathcal{B}(S)$ is its Borel $\sigma$-algebra.

## 2 Fock Space Preliminaries

In this section we define some key operators in Fock space. Full proofs of all results mentioned here can be found in [43], [38] or [34].

Let $H$ be a complex separable Hilbert space and $H^{\otimes n}$ be its $n$-fold tensor product. We denote by $H_s^{\otimes n}$ the closed subspace of $H^{\otimes n}$ comprising symmetric tensors. $H_s^{\odot n}$ is the image of $H_s^{\otimes n}$ under the bijection $\psi \to \sqrt{n!}\psi$ and is regarded as a Hilbert space with respect to the inner product $\langle \cdot, \cdot \rangle_{H_s^{\odot n}} = n!\langle \cdot, \cdot \rangle_{H_s^{\otimes n}}$. We write $f^{\odot n} := \sqrt{n!} f^{\otimes n}$. If $g \in H$ we define the symmetrisation of $f^{\odot n}$ and $g$ to be the vector $\text{Symm}(f^{\odot n}, g) \in H_s^{\odot n+1}$ defined by

$$\text{Symm}(f^{\odot n}, g) := \sqrt{\frac{n!}{n+1}} \sum_{r=0}^{n} f^{\otimes n-r} \otimes g \otimes f^{\otimes r}.$$

Note that the choice of normalisation ensures that $\text{Symm}(f^{\odot n}, f) = f^{\odot n+1}$.

Symmetric Fock space over $H$ is $\Gamma(H) := \bigoplus_{n=0}^{\infty} H_s^{\otimes n}$, where by convention $H_s^{\otimes 0} := \mathbb{C}$ and $H_s^{\otimes 1} := H$. We also define $\hat{\Gamma}(H) := \bigoplus_{n=0}^{\infty} H_s^{\odot n}$. Of course $\hat{\Gamma}(H)$ and $\Gamma(H)$ are naturally isomorphic (see below) and we will freely move between these spaces in the sequel. We will often identify $H_s^{\otimes n}$ with its natural embedding in $\Gamma(H)$ whereby each $f_n \in H_s^{\otimes n}$ is mapped to $(0, \ldots, 0, f_n, 0, \ldots)$.

The exponential vector $e(f) \in \Gamma(H)$ corresponding to $f \in H$ is defined by $e(f) = \left(1, f, \frac{f \otimes f}{\sqrt{2!}}, \ldots, \frac{f^{\otimes n}}{\sqrt{n!}}, \ldots\right)$. We have $\langle e(f), e(g) \rangle = e^{\langle f, g \rangle}$, for each $f, g \in H$. The mapping $f \to e(f)$ from $H$ to $\Gamma(H)$ is analytic.

If $\mathcal{D} \subseteq H$ we define $\mathcal{E}(\mathcal{D})$ to be the linear span of $\{e(f), f \in \mathcal{D}\}$. In particular if $\mathcal{D}$ is a dense linear manifold in $H$ then $\mathcal{E}(\mathcal{D})$ is dense in $\Gamma(H)$. We define $\mathcal{E} := \mathcal{E}(H)$. If $f \in H$, the corresponding annihilation operator $a(f)$, creation operator $a^\dagger(f)$, exponential annihilation operator $U(f)$ and exponential creation operator $U^\dagger(f)$ are defined on $\mathcal{E}$ by linear extension of the following prescriptions:

$$a(f)e(g) = \langle f, g \rangle e(g), \tag{2.1}$$

$$a^\dagger(f)e(g) = \left.\frac{d}{dt}e(g + tf)\right|_{t=0}, \tag{2.2}$$

$$U(f)e(g) = e^{\langle f,g \rangle}e(g), \tag{2.3}$$

$$U^\dagger(f)e(g) = e(g + f), \tag{2.4}$$

for all $g \in H$. Each of these operators is closable, indeed $a^\dagger(f) \subseteq a(f)^*, U^\dagger(f) \subseteq U(f)^*$, for each $f \in H$. We also have the canonical commutation relations

$$a(f)a^\dagger(g)\psi - a^\dagger(g)a(f)\psi = \langle f, g \rangle \psi, \tag{2.5}$$

for all $f, g \in H, \psi \in \mathcal{E}$.

If $T$ is a contraction in $H$ then its second quantisation $\Gamma(T)$ is the contraction in $\Gamma(H)$ whose action on $\mathcal{E}$ is given by linear extension of

$$\Gamma(T)e(f) = e(Tf). \tag{2.6}$$

In particular if $T$ is unitary, then so is $\Gamma(T)$. If $A$ is a bounded self-adjoint operator in $H$ we define the associated conservation operator $\Lambda(A)$ to be the infinitesimal generator of the one-parameter unitary group $(\Gamma(e^{itA}), t \in \mathbb{R})$. It is easily checked that $\mathcal{E} \subseteq \mathrm{Dom}(\Lambda(A))$ and we have the useful identity

$$\langle e(f), \Lambda(A)e(g) \rangle = \langle f, Ag \rangle \langle e(f), e(g) \rangle, \tag{2.7}$$

for all $f, g \in H$.

The number operator is defined by $N := \Lambda(I)$. Its domain is $\mathrm{Dom}(N) = \{(f_n, n \in \mathbb{Z}_+) \in \Gamma(H); \sum_{n=1}^\infty n^2 ||f_n||^2 < \infty\}$. The associated contraction semigroup is $(T_t, t \geq 0)$ where for each $t \geq 0$

$$T_t := e^{-tN} = \sum_{n=0}^\infty e^{-tn} P_n,$$

and where $P_n$ denotes the orthogonal projection from $\Gamma(H)$ to $H_s^{\otimes n}$.

The unitary isomorphism from $\Gamma(H)$ to $\hat{\Gamma}(H)$ which we employ to identify these spaces is $(N!)^{-\frac{1}{2}}$.

Using the analyticity of exponential vectors, we can easily check that finite particle vectors are in the domains of annihilation and creation operators and we deduce the following:

$$a(f)g^{\odot n} = n\langle f, g\rangle g^{\odot n-1}, \qquad (2.8)$$

$$a^\dagger(f)g^{\odot n} = \mathrm{Symm}(g^{\odot n}, f), \qquad (2.9)$$

for all $f, g \in H$.

If $H = H_1 \oplus H_2$ we may identify $\Gamma(H)$ with $\Gamma(H_1) \otimes \Gamma(H_2)$ via the natural isomorphism which maps $e(f)$ to $e(f_1) \otimes e(f_2)$ for each $f = (f_1, f_2) \in H$. In this context we will always denote the linear span of the exponential vectors in $H_i$ by $\mathcal{E}_i (i = 1, 2)$.

## 3 Universal Annihilation and Creation Operators

For each $t \in \mathbb{R}$ we define linear operators $V_t : \Gamma(H) \to \Gamma(H) \otimes \Gamma(H)$ and $V_t^\dagger : \Gamma(H) \otimes \Gamma(H) \to \Gamma(H)$ on the dense domains $\mathcal{E}$ and $\mathcal{E} \otimes \mathcal{E}$ (respectively) by linear extension of the following prescriptions:

$$V_t e(f) = e(f) \otimes e(tf), \qquad (3.10)$$

$$V_t^\dagger(e(f) \otimes e(g)) = U(tg)^\dagger e(f) = e(f + tg), \qquad (3.11)$$

for all $f, g \in H$.

**Proposition 3.1** *For each $t \in \mathbb{R}$, $V_t$ and $V_t^\dagger$ are closable with $V_t^\dagger \subseteq V_t^*$.*

*Proof.* The result will follow if we can show that these operators are mutually adjoint. Now for all $f, g, h \in H, t \in \mathbb{R}$,

$$\begin{aligned}
\langle V_t e(f), e(g) \otimes e(h)\rangle &= \langle e(f) \otimes e(tf), e(g) \otimes e(h)\rangle \\
&= \langle e(f), e(g)\rangle \langle e(tf), e(h)\rangle \\
&= \langle e^{t\langle h, f\rangle} e(f), e(g)\rangle \\
&= \langle U(th)e(f), e(g)\rangle \quad \text{(by 2.3)} \\
&= \langle e(f), U(th)^\dagger e(g)\rangle \\
&= \langle e(f), V_t^\dagger e(g) \otimes e(h)\rangle.
\end{aligned}$$

□

We define linear operators $\nabla^- : \Gamma(H) \to \Gamma(H) \otimes H$ and $\nabla^+ : \Gamma(H) \otimes H \to \Gamma(H)$ on the dense domains $\mathcal{E}$ and $\mathcal{E} \otimes H$ (respectively) by linear extension of the following prescriptions, for each $f, g \in H$:

$$\nabla^- e(f) = \frac{d}{dt} V_t e(f) \Big|_{t=0},$$

so by (3.10)
$$\nabla^- e(f) = e(f) \otimes f, \qquad (3.12)$$

and $\quad \nabla^+(e(f) \otimes g) = \frac{d}{dt} V_t^\dagger e(f) \otimes e(g) \Big|_{t=0},$

so by (3.11)
$$\nabla^+(e(f) \otimes g) = a^\dagger(g) e(f). \qquad (3.13)$$

**Proposition 3.2** $\nabla^-$ and $\nabla^+$ are closable with $\nabla^+ \subseteq (\nabla^-)^*$.

*Proof.* The fact that $\nabla^-$ and $\nabla^+$ are mutually adjoint follows from differentiation of the adjunction relation between $V_t$ and $V_t^\dagger$ and the result follows. $\square$

We call $\nabla^-$ and $\nabla^+$ *universal annihilation* and *universal creation* operators (respectively). For each $n \in \mathbb{Z}_+$ we denote the restrictions of $\nabla^-$ and $\nabla^+$ to $H_s^{\odot n}$ and $H_s^{\odot n} \otimes H$ (respectively) by $\nabla_n^-$ and $\nabla_n^+$. Using the analyticity of exponential vectors, we can easily deduce that $\text{Ran}(\nabla_n^-) \subseteq H_s^{\odot n} \otimes H$ and $\text{Ran}(\nabla_n^+) \subseteq H_s^{\odot n+1}$ and obtain the following analogues of (2.8) and (2.9) (c.f. [47]):

$$\nabla_n^- f^{\odot n} = n f^{\odot n-1} \otimes f \qquad (3.14)$$

$$\nabla_n^+ f^{\odot n} \otimes g = \text{Symm}(f^{\odot n}, g), \qquad (3.15)$$

for each $f, g \in H$.

**Lemma 3.1** For each $n \in \mathbb{Z}_+$, $\nabla_n^-$ and $\nabla_n^+$ are bounded operators with $||\nabla_n^-|| = \sqrt{n}$ and $||\nabla_n^+|| = \sqrt{n+1}$.

*Proof.* For each $f \in H$,
$$\begin{aligned} ||\nabla_n^- f^{\odot n}|| &= n || f^{\odot n-1} \otimes f|| \\ &= n\sqrt{(n-1)!} ||f||^n = \sqrt{n} ||f^{\odot n}||_{H_s^{\odot n}}. \end{aligned}$$

The result for $\nabla_n^-$ follows from the fact that the linear span of $\{f^{\odot n}, f \in H\}$ is dense in $H_s^{\odot n}$.

For each $f, g \in H$,

$$\begin{aligned}||\nabla_n^+(f^{\odot n} \otimes g)|| &= \left\|\sqrt{\frac{n!}{n+1}}\sum_{r=0}^n f^{\otimes n-r} \otimes g \otimes f^{\otimes r}\right\| \\ &\leq \sqrt{(n+1)!}||f||^n||g|| \\ &= \sqrt{n+1}||f^{\odot n} \otimes g||_{H_s^{\odot n} \otimes H}.\end{aligned}$$

Since the linear span of $\{f^{\odot n} \otimes g, f, g \in H\}$ is dense in $H_s^{\odot n} \otimes H$, it follows that $\nabla_n^+$ is bounded. To see that the bound is obtained, observe that for all $f \in H$,

$$\nabla_n^+(f^{\odot n} \otimes f) = \text{Symm}(f^{\odot n}, f) = \sqrt{n+1} f^{\odot n} \otimes f. \qquad \square$$

If $\phi_n \in H_s^{\odot n} \otimes H$, there exist sequences $(f_{n,r}, r \in \mathbb{N})$ and $(g_r, r \in \mathbb{N})$ where each $f_{n,r} \in H_s^{\odot n}$ and $g_r \in H$ such that $\phi_n = \sum_{r=1}^\infty f_{n,r} \otimes g_r$. We define

$$\widetilde{\phi_n} := \nabla_n^+ \phi_n = \sum_{r=1}^\infty \text{Symm}(f_{n,r}, g_r).$$

In the next theorem we will find it convenient to identify $\hat{\Gamma}(H) \otimes H$ with $\bigoplus_{n=0}^\infty (H_s^{\odot n} \otimes H)$.

**Theorem 3.1** *1. $\text{Dom}(\nabla^-) = \left\{\psi = (\psi_n, n \in \mathbb{Z}_+) \in \hat{\Gamma}(H); \sum_{n=1}^\infty nn!||\psi_n||^2 < \infty\right\}$.*

*2. $\text{Dom}(\nabla^+) = \left\{\phi = (\phi_n, n \in \mathbb{Z}_+) \in \hat{\Gamma}(H) \otimes H; \sum_{n=0}^\infty ||\widetilde{\phi_n}||^2 < \infty\right\}$.*

*Proof.*

1. (*Sufficiency*)

   Let $\psi \in \hat{\Gamma}(H)$ be such that $\sum_{n=1}^\infty nn!||\psi_n||^2 < \infty$ and for each $M \in \mathbb{Z}_+$, define $\psi^{(M)} = (\psi_0, \psi_1, \ldots, \psi_M, 0, 0, \ldots)$. Clearly $\psi^{(M)} \to \psi$ as $M \to \infty$. Using Lemma 3.1 we see that for each $M, N \in \mathbb{N}, N > M$

   $$||\nabla^-\psi_N - \nabla^-\psi_M||^2 = \sum_{n=M+1}^N nn!||\psi_n||^2 \to 0 \text{ as } N, M \to \infty.$$

   Hence $(\nabla^-\psi_N, N \in \mathbb{N})$ converges to a vector $\phi$ in $\hat{\Gamma}(H) \otimes H$. Since $\nabla^-$ is closed, we deduce that $\phi = \nabla^-\psi$ and so $\psi \in \text{Dom}(\nabla^-)$.

*(Necessity)* Suppose that $\psi = (\psi_n, n \in \mathbb{Z}_+) \in \text{Dom}(\nabla^-)$, then $\nabla^- \psi = (\nabla_n^- \psi_n, n \in \mathbb{Z}_+)$ and again using Lemma 3.1 we obtain

$$\begin{aligned} ||\nabla^- \psi||^2 &= \sum_{n=0}^{\infty} ||\nabla_n^- \psi_n||^2 \\ &= \sum_{n=1}^{\infty} nn! ||\psi_n||^2 < \infty. \end{aligned}$$

2. This is proved by exactly the same argument as (1). $\square$

Clearly, as an operator on $\hat{\Gamma}(H)$, $\text{Dom}(\nabla^-) = \text{Dom}(\sqrt{N})$ (see [34], [33]). To some extent, Theorem 3.1 (2) tells us less than (1) (although it is sufficient for applications in probability). Another approach to $\nabla^+$ is given by J.M.Lindsay in [34]. Let $\Phi(H)$ be the full Fock space over $H$. Then since $\Gamma(H) \otimes H = \bigoplus_{n=0}^{\infty} (H_s^{\otimes n} \otimes H) \subseteq \Phi(H)$ we can regard $\nabla^+$ as an operator from $\Phi(H)$ to $\Gamma(H)$. In fact it is not difficult to verify that in this case $\nabla^+ = \sqrt{N} P_s$ where $P_s$ is the orthogonal projection from $\Phi(H)$ to $\Gamma(H)$.

Now suppose that $H = H_1 \oplus H_2$, then as previously remarked we may identify $\Gamma(H)$ with $\Gamma(H_1) \otimes \Gamma(H_2)$. We may also identify $\Gamma(H) \otimes H$ with $[(\Gamma(H_1) \otimes H_1) \otimes \Gamma(H_2)] \oplus [\Gamma(H_1) \otimes (\Gamma(H_2) \otimes H_2)]$ in an obvious way. For $i = 1, 2$ $\nabla_i^{\pm}$ denotes the universal annihilation/creation operators associated to each $H_i$ and $\pi_i$ are the isometric embeddings of $H_i$ into $H_1 \oplus H_2$, so for example if $f \in H_1$ then $\pi(f) = (f, 0)$. Note that $\pi_1^*((f_1, f_2)) = f_1$ for all $f_i \in H_i$. For simplicity, we will continue to use the notation $\pi_i$ and $\pi_i^*$ when these operators are tensored with the identity to act in tensor products.

Using the identifications given above, it is not difficult to verify that

$$\nabla^- = \pi_1(\nabla_1^- \otimes I) + \pi_2(I \otimes \nabla_2^-), \qquad (3.16)$$

on $\text{Dom}(\nabla_1^-)\underline{\otimes}\text{Dom}(\nabla_2^-)$, and

$$\nabla^+ = (\nabla_1^+ \otimes I)\pi_1^* + (I \otimes \nabla_2^+)\pi_2^*, \qquad (3.17)$$

on $\text{Dom}(\nabla_1^+) \oplus \text{Dom}(\nabla_2^+)$ (c.f. [43], p.150).

It is useful to think of $\nabla^-$ as a "gradient" and $\nabla^+$ as a "divergence" and we will make these correspondences precise later. In this respect, we should define associated "directional derivatives". For this purpose we introduce the Dirac "bra" and "ket" operators $\epsilon_f : H \to \mathbb{C}$ and $\epsilon_f^\dagger : \mathbb{C} \to H$ by

$$\epsilon_f(g) = \langle f, g \rangle, \quad \epsilon_f^\dagger(\alpha) = \alpha f,$$

for each $f, g \in H, \alpha \in \mathbb{C}$ (c.f. [16]). These operators are clearly linear, bounded and mutually adjoint with each $||\epsilon_f|| = ||\epsilon_f^\dagger|| = ||f||$. The nature of the "directional derivative" operators is revealed in the following result:

**Proposition 3.3** *For each $f \in H$,*

1. $(I \otimes \epsilon_f) \circ \nabla^- = a(f)$,

2. $\nabla^+ \circ (I \otimes \epsilon_f^\dagger) = a^\dagger(f)$, [1]

*on $\mathcal{E}(H)$.*

Proof.

1. For each $f, g \in H$,
$$(I \otimes \epsilon_f) \circ \nabla^- e(g) = (I \otimes \epsilon_f)(e(g) \otimes g) = \langle f, g \rangle e(g) = a(f) e(g).$$

2. follows by taking adjoints in (1).

$\square$

Using a density argument, it is easily verified that Proposition 3.3 (1) extends to $\mathrm{Dom}(\nabla^-)$. Moreover it follows that $\mathrm{Dom}(\nabla^-)$ is the maximal domain for all $a(f), f \in H$.

We will now give a noncommutative factorisation of the number operator (see [47]) for a different factorisation in the additive sense). First, for each $n \in \mathbb{N}$, we define a linear operator $W_n$ from $H^{\odot n}$ to $H^{\odot n-1} \otimes H$ by
$$W_n = \frac{1}{\sqrt{n}} \nabla_n^-.$$

We also define $W : \hat{\Gamma}(H) \to \hat{\Gamma}(H) \otimes H$ by $W = \bigoplus_{n=0}^{\infty} W_n$ where $W_0 := 1$.

**Theorem 3.2**  1. *$W_n$ is unitary for each $n \in \mathbb{Z}_+$.*

2. *$W$ is unitary.*

3. *On their maximal domains,*
$$\nabla^- = W\sqrt{N} = (\sqrt{N+1} \otimes I)W,$$
$$\nabla^+ = \sqrt{N}W^* = W^*(\sqrt{N+1} \otimes I).$$

---

[1] Here we are identifying $H$ with $H \otimes \mathbb{C}$.

4. On $Dom(N)$,
$$N = \nabla^+ \nabla^-.$$

*Proof.*

1. Let $f \in H$ then $W_n f^{\odot n} = \sqrt{n} f^{\odot n-1} \otimes f$ and so
$$\begin{aligned}||W_n f^{\odot n}|| &= \sqrt{n} ||f^{\odot n-1}|| \cdot ||f|| \\ &= \sqrt{n!} ||f||^n = ||f^{\odot n}||.\end{aligned}$$

   Hence $W_n$ is an isometry between total sets and so extends to a unitary operator by linearity and continuity.

2. Follows immediately from (1).

3. $\nabla^- = W\sqrt{N}$ on $Dom(\nabla^-) = Dom(\sqrt{N})$ is immediate. Since $W_n^* = \frac{1}{\sqrt{n}} \nabla_n^+$ it follows that $W^* = N^{-\frac{1}{2}} \nabla^+$ and hence $W^*(Dom(\nabla^+)) \subseteq Dom(\sqrt{N})$. So $\nabla^+ = \sqrt{N} W^*$ on $Dom(\nabla^+)$. The other results are proved similarly.

4. This follows immediately from (3) and (2). $\square$

$Dom(\nabla^-)$ becomes a complex Hilbert space with respect to the inner product $\langle \cdot, \cdot \rangle_1$ where for each $\psi_1, \psi_2 \in Dom(\nabla^-)$,
$$\langle \psi_1, \psi_2 \rangle_1 = \langle \psi_1, \psi_2 \rangle + \langle \nabla^- \psi_1, \nabla^- \psi_2 \rangle.$$

Consider the self-adjoint linear operator $Q := (1+N)^{-\frac{1}{2}} = \pi^{-\frac{1}{2}} \int_0^\infty t^{-\frac{1}{2}} e^{-t} T_t dt$, (see Lemma 3.12 in [23] for the last identity).

It is easy to check that $Q$ is a bounded operator on $\Gamma(H)$ and that $Q Dom(\nabla^-) \subseteq Dom(N)$.

**Theorem 3.3** *$Q$ is a unitary isomorphism between $(\Gamma(H), ||\cdot||)$ and $(Dom(\nabla^-), ||\cdot||_1)$.*

*Proof.* (c.f. the proof of Proposition 3.14 in [23]). Let $\psi_1, \psi_2 \in Dom(\nabla^-)$, then by Theorem 3.2
$$\begin{aligned}\langle \psi_1, \psi_2 \rangle &= \langle Q^{-1} Q \psi_1, Q^{-1} Q \psi_2 \rangle \\ &= \langle Q^{-2} Q \psi_1, Q \psi_2 \rangle \\ &= \langle Q \psi_1, Q \psi_2 \rangle + \langle N Q \psi_1, Q \psi_2 \rangle \\ &= \langle Q \psi_1, Q \psi_2 \rangle + \langle \nabla^- Q \psi_1, \nabla^- Q \psi_2 \rangle \\ &= \langle Q \psi_1, Q \psi_2 \rangle_1.\end{aligned}$$

Hence $Q$ is an isometric embedding of $(\Gamma(H), ||\cdot||)$ in $(\text{Dom}(\nabla^-), ||\cdot||_1)$. The result follows from the fact that $\text{Ran}(Q)$ is dense in $(\text{Dom}(\nabla^-), ||\cdot||_1)$. To see this let $\mathcal{D}$ be the linear space comprising those sequences $(\psi_n, n \in \mathbb{N}) \in \Gamma(H)$ where $\psi_n = 0$ for all but finitely many $n$. Clearly $\mathcal{D}$ is dense in $\text{Dom}(\nabla^-)$. However $Q\mathcal{D} = \mathcal{D}$ and the result follows. $\square$

From now on we will use $\Xi$ to denote the Hilbert space $\text{Dom}(\nabla^-)$ equipped with the inner product $\langle \cdot, \cdot \rangle_1$. In the sequel we will also want to work with the closed linear operator $\nabla^- \otimes I$ acting in $\Gamma(H) \otimes H$. The Hilbert space $\text{Dom}(\nabla^- \otimes I)$ equipped with the graph norm is precisely $\Xi \otimes H$.

Let $\tau$ be the tensor shift on $H \otimes H$ so that $\tau$ is the closed linear extension of the map $\tau(f \otimes g) = g \otimes f$ for each $f, g \in H$. It is easily verified that $\tau$ is self-adjoint and unitary (in physicists' language, $\tau$ is an example of a "parity operator").

Our next result is an abstract version of the "Itô-Skorohod isometry". Note however that the presence of the operator $\tau$ ensures that the isometry property does not in fact hold.

**Theorem 3.4** *For all $\psi_i \in \Xi \otimes H (i = 1, 2)$*

$$\langle \nabla^+ \psi_1, \nabla^+ \psi_2 \rangle = \langle \psi_1, \psi_2 \rangle + \langle (I \otimes \tau)(\nabla^- \otimes I)\psi_1, (\nabla^- \otimes I)\psi_2 \rangle. \quad (3.18)$$

*Furthermore $\nabla^+$ is a contraction from $\Xi \otimes H$ into $\Gamma(H)$.*

*Proof.* To establish (3.18), let $f_i, g_i \in H (i = 1, 2)$. We then find that for each $n \in \mathbb{N}$,

$$\begin{aligned}
\langle \nabla_n^+ f_1^{\odot n} \otimes g_1, \nabla_n^+ f_2^{\odot n} \otimes g_2 \rangle &= \frac{n!}{n+1} \sum_{r=0}^{n} \sum_{s=0}^{n} \langle f_1^{\otimes n-r} \otimes g_1 \otimes f_1^{\otimes r}, f_2^{\otimes n-s} \otimes g_2 \otimes f_2^{\otimes s} \rangle \\
&= \langle f_1^{\odot n} \otimes g_1, f_2^{\odot n} \otimes g_2 \rangle + n^2 \langle f_1^{\odot n-1}, f_2^{\odot n-1} \rangle \langle g_1, f_2 \rangle \langle f_1, g_2 \rangle \\
&= \langle f_1^{\odot n} \otimes g_1, f_2^{\odot n} \otimes g_2 \rangle \\
&\quad + \langle (I \otimes \tau)(\nabla_n^- \otimes I) f_1^{\odot n} \otimes g_1, (\nabla_n^- \otimes I) f_2^{\odot n} \otimes g_2 \rangle.
\end{aligned}$$

The required result follows from here by linearity and continuity. To establish the contraction property, we have from (3.18) that for all $\psi \in \Xi \otimes H$,

$$||\nabla^+ \psi||^2 = ||\psi||^2 + \langle (I \otimes \tau)(\nabla^- \otimes I)\psi, (\nabla^- \otimes I)\psi \rangle,$$

and the result follows easily from this by using the Cauchy-Schwarz inequality and the isometry property of $I \otimes \tau$. $\square$

The last result is closely related to the canonical commutation relations (2.5). This is not so clear from the argument in the proof of Theorem 3.4 however it is instructive to compute actions on exponential vectors. If $f_i, g_i \in H (i = 1, 2)$, we find that

$$\begin{aligned}
\langle \nabla^+ e(f_1) \otimes g_1, \nabla^+ e(f_2) \otimes g_2 \rangle &= \langle a^\dagger(g_1)e(f_1), a^\dagger(g_2)e(f_2) \rangle \\
&= \langle e(f_1), e(f_2) \rangle \langle g_1, g_2 \rangle + \langle a(g_2)e(f_1), a(g_1)e(f_2) \rangle \\
&= \langle e(f_1) \otimes g_1, e(f_2) \otimes g_2 \rangle \\
&+ \langle (I \otimes \tau)(\nabla^- \otimes I)e(f_1) \otimes g_1, (\nabla^- \otimes I)e(f_2) \otimes g_2 \rangle.
\end{aligned}$$

Now let $H = L^2(S, \mathcal{S}, \mu)$ where $S$ is a locally compact topological space, $\mathcal{S}$ is its Borel $\sigma$-algebra and $\mu$ is a Borel measure defined on $(S, \mathcal{S})$. In this case $\Gamma(H) \otimes H = L^2(S, \mathcal{S}, \mu; \Gamma(H))$. If $\psi_n \in H^{\odot n}$, $\psi_{n-1}(\cdot, s)$ will denote the symmetric function of $n - 1$ variables obtained by fixing $s \in S$. We then have for all $\psi = (\psi_n, n \in \mathbb{Z}_+) \in \text{Dom}(\nabla^-)$,

$$||\nabla^- \psi||^2 = \int_S ||\nabla_s^- \psi||^2_{\Gamma(H)} \mu(ds),$$

where for $\mu$-almost all $s \in S$

$$\nabla_s^- \psi := (n \psi_{n-1}(\cdot, s), n \in \mathbb{N}). \tag{3.19}$$

We call $\nabla_s^- \psi$ the *Lindsay-Malliavin transform* of $\psi$ at $s$. Note that $\nabla_s^-$ is not a bona fide operator since the right hand side of (3.19) depends on the choice of a sequence of functions from the equivalence class of $\psi$.

If $\mu$ is a regular measure (so that compact sets have finite mass) then the space $\mathcal{D}$ of continuous functions with compact support is a dense subspace of $H$. In this case the Lindsay-Malliavin transform is a genuine linear operator whose action on $\mathcal{E}(\mathcal{D})$ is given by

$$\nabla_s^- e(f) = f(s)e(f),$$

for all $f \in \mathcal{D}$. So each $\nabla_s^-$ is densely defined on $\mathcal{E}(\mathcal{D})$; however it is well known that these operators are not closable (see e.g. [32]).

Even when $\mu$ fails to be regular we can rewrite the result of Theorem 3.4 by using the Lindsay-Malliavin transform. First we observe that $\Xi \otimes H = L^2(S, \mathcal{S}, \mu; \Xi)$ and elements of this space may be regarded as equivalence classes of mappings from $S$ to $\Xi$. Note that $\nabla^- \otimes I : L^2(S, \mathcal{S}, \mu; \Xi) \to L^2(S^2, \mathcal{S}^{\otimes 2}, \mu \times \mu; \Gamma(H))$.

**Corollary 3.1** *If $X, Y \in L^2(S, \mathcal{S}, \mu; \Xi)$, then*

$$\langle \nabla^+ X, \nabla^+ Y \rangle = \int_S \langle X(s), Y(s) \rangle \mu(ds) + \int_S \int_S \langle \nabla_t^- X(s), \nabla_s^- Y(t) \rangle \mu(ds) \mu(dt). \tag{3.20}$$

*Proof* We use the same notation as in the proof of Theorem 3.4. It is sufficient to consider the case where $X(t) = f_1^{\odot^n} g_1(t)$ and $Y(t) = f_2^{\odot^n} g_2(t)$ for each $t \in S$. Using (3.19) we obtain

$$\begin{aligned}
& \langle (I \otimes \tau)(\nabla_n^- \otimes I) f_1^{\odot^n} \otimes g_1, (\nabla_n^- \otimes I) f_2^{\odot^n} \otimes g_2 \rangle \\
&= n^2 \langle f_1^{\odot^{n-1}}, f_2^{\odot^{n-1}} \rangle \langle g_1, f_2 \rangle \langle f_1, g_2 \rangle \\
&= n^2 \langle f_1^{\odot^{n-1}}, f_2^{\odot^{n-1}} \rangle \int_S \overline{g_1(s)} f_2(s) \mu(ds) \int_S \overline{f_1(t)} g_2(t) \mu(ds) \\
&= \int_S \int_S \langle \nabla_t^- X(s), \nabla_s^- Y(t) \rangle \mu(ds) \mu(dt),
\end{aligned}$$

and the result follows. $\square$

**Remark.** Some of the main results of this section - Theorems 3.1, 3.2 (4), 3.4 and Corollary 3.1 are all given, at least in outline in [34] Propositions 3.1 and 3.2 (see also [33] for the Guichardet space version). A similar result to Theorem 3.4 is also established in [48] - see Proposition 1 therein.

# 4 The Chaos Decomposition of Lévy-Itô Space

## 4.1 Preliminaries on Lévy Processes [2]

Let $(\Omega, \mathcal{F}, (\mathcal{F}_t, t \geq 0), P)$ be a stochastic base wherein the filtration $(\mathcal{F}_t, t \geq 0)$ satisfies the usual hypotheses of completeness and right continuity. Let $X = (X(t), t \geq 0)$ be an adapted real-valued Lévy process defined on $(\Omega, \mathcal{F}, P)$ so that $X(0) = 0$ (a.s.), $X$ has stationary increments and strongly independent increments (in the sense that $X(t) - X(s)$ is independent of $\mathcal{F}_s$ for all $0 \leq s < t < \infty$), $X$ is stochastically continuous and its paths are a.s. càdlàg. We have the Lévy-Khintchine formula

$$\mathbb{E}(e^{iuX(t)}) = e^{-t\eta(u)},$$

for all $t \geq 0, u \in \mathbb{R}$, where $\eta : \mathbb{R} \to \mathbb{C}$ is a continuous, hermitian negative definite mapping for which $\eta(0) = 0$. It has the canonical form

$$\begin{aligned}
\eta(u) &= -ibu + \frac{1}{2}\sigma^2 u^2 \\
&+ \int_{\mathbb{R}-\{0\}} (1 - e^{iuy} + iuy1_{B_1}(y))\nu(dy),
\end{aligned}$$

where $b \in \mathbb{R}, \sigma \geq 0$ and $\nu$ is a Lévy measure on $\mathbb{R} - \{0\}$, i.e. $\nu$ is a Borel measure for which $\int_{\mathbb{R}-\{0\}} (1 \wedge |y|^2)\nu(dy) < \infty$. Information about the sample paths of $X$ is given by the *Lévy-Itô decomposition*:

$$X(t) = bt + \sigma B(t) + \int_{|x|<1} x\tilde{N}(t,dx) + \int_{|x|\geq 1} xN(t,dx). \quad (4.21)$$

Here $N$ is the Poisson random measure on $\mathbb{R}^+ \times (\mathbb{R} - \{0\})$ defined by

$$N(t,A) := \#\{0 \leq s \leq t, \Delta X(s) \in A\},$$

for each $t \geq 0, A \in \mathcal{B}(\mathbb{R} - \{0\})$, $\tilde{N}$ is the compensated random measure defined by

$$\tilde{N}(t,A) = N(t,A) - t\nu(A),$$

and $B = (B(t), t \geq 0)$ is a standard Brownian motion which is independent of $N$.

## 4.2 The Itô Representation Theorem

In this and the next section, we briefly outline proofs of results which are given more fully in [2](ii).

We fix $T > 0$ and let $f \in L^2([0,T], \mathbb{R})$. We may then form the *Wiener-Itô integral* $X_f(t) = \int_0^t f(s)dX(s)$, for each $0 \leq t \leq T$. We define

$$M_f(t) := \exp\left\{iX_f(t) + \int_0^t \eta(f(s))ds\right\}.$$

**Lemma 4.1** *For each $f \in L^2([0,T]), u \in \mathbb{R}, t \in [0,T]$,*

1. $\mathbb{E}(e^{iuX_f(t)}) = \exp\left\{-\int_0^t \eta(uf(s))ds\right\}.$

2. *$(M_f(t), t \in [0,T])$ is a complex-valued square-integrable martingale with stochastic differential*

$$dM_f(t) = i\sigma f(t)M_f(t-)dB(t) + (e^{if(t)x} - 1)M_f(t-)\tilde{N}(dt,dx). \quad (4.22)$$

*Proof.* These are both straightforward applications of Itô's formula applied to the processes $(e^{iuX_f(t)}, 0 \leq t \leq T)$ and $(M_f(t), 0 \leq t \leq T)$, respectively. $\square$

From now on, for each $0 \leq t \leq T$, we require that $\mathcal{F}_t = \sigma\{X(s), 0 \leq s \leq t\}$. We define *Lévy-Itô space* to be the complex separable Hilbert space $\mathcal{H} := L^2(\Omega, \mathcal{F}_T, P; \mathbb{C})$. $\mathcal{P}_T$ will denote the predictable $\sigma$-algebra generated by processes defined on $[0,T] \times \Omega$.

**Lemma 4.2** $\{M_f(T), f \in L^2([0,T])\}$ *is total in* $\mathcal{H}$.

*Proof.* This is a consequence of the injectivity of the Fourier transform. The proof is similar to that of Lemma 4.3.2 in [42]. □

Let $\mathcal{H}_2^{(B)}(T)$ be the complex Hilbert space of all complex predictable processes satisfying $\int_0^T \mathbb{E}(|F(t)|^2)dt < \infty$ and $\mathcal{H}_2^{(N)}(T)$ be the complex Hilbert space of all $\mathcal{P}_T \times \mathcal{B}(\mathbb{R}-\{0\})$ measurable mappings $G: [0,T] \times (\mathbb{R}-\{0\}) \times \Omega \to \mathbb{C}$ for which $\int_0^T \int_{\mathbb{R}-\{0\}} \mathbb{E}(|G(t,x)|^2)\nu(dx)dt < \infty$.

**Theorem 4.1** *[The Itô Representation]*

If $F \in \mathcal{H}$, then there exists unique $\psi_0 \in \mathcal{H}_2^{(B)}(T)$ and $\psi_1 \in \mathcal{H}_2^{(N)}(T)$ such that

$$F = \mathbb{E}(F) + \sigma \int_0^T \psi_0(s) dB(s) + \int_0^T \int_{\mathbb{R}-\{0\}} \psi_1(s,x) \tilde{N}(ds,dx). \qquad (4.23)$$

*Proof.* The result holds for $F = M_f(T)$ by (4.22) and is easily extended to finite linear combinations of such random variables. The extension to arbitrary $F$ is by approximation using Lemma 4.2 (see the proof of Theorem 4.3.3. in [42], Proposition 3 in [35]) and Proposition 2.1 in [44].) □

## 4.3 Multiple Wiener-Lévy Integrals and the Chaos Decomposition

Let $X$ be a Lévy process with associated Lévy-Itô decomposition (4.21). Let $S = [0,T] \times \mathbb{R}$. We consider the associated martingale-valued measure $M$ defined on $(S, \mathcal{I})$ by the prescription

$$M([0,t] \times A) = \tilde{N}(t, A-\{0\}) + \sigma B(t)\delta_0(A)$$

for each $t \in [0,T], A \in \mathcal{B}(\mathbb{R})$ where $\mathcal{I}$ is the ring comprising finite unions of sets of the form $I \times A$ where $A \in \mathcal{B}(\mathbb{R})$ and $I$ is itself a finite union of intervals (see e.g. [2] for more information about martingale-valued measures). The associated "control measure" is the $\sigma$-finite measure $\mu = \lambda \times \rho$, where $\lambda$ is Lebesgue measure on $[0,T]$ and $\rho$ is defined on $(\mathbb{R}, \mathcal{B}(\mathbb{R}))$ by $\rho(A) = \sigma^2 \delta_0(A) + \nu(A-\{0\})$ for all $A \in \mathcal{B}(\mathbb{R})$. We can easily compute

$$\mathbb{E}(M([0,t] \times A)^2) = \mu([0,t] \times A) = t\rho(A),$$

for each $t \in [0,T], A \in \mathcal{B}(\mathbb{R})$.

Returning to the set-up of section 1, we take $H = L^2(S, \mathcal{B}(S), \mu; \mathbb{C})$ so that for each $n \in \mathbb{N}$, $H^{\otimes n} = L^2(S^n, \mathcal{B}(S^n), \mu^n; \mathbb{C})$ and $H_s^{\otimes n}$ comprises symmetric complex-valued square-integrable functions on $S^n$. Fix $n \in \mathbb{N}$ and define $\mathcal{D}^{(n)}$ to be the linear space of all functions $f_n \in H^{\otimes n}$ which take the form

$$f_n = \sum_{j_1,\ldots,j_n=1}^{N} a_{j_1,\ldots,j_n} 1_{A_{j_1} \times \cdots \times A_{j_n}}, \quad (4.24)$$

where $N \in \mathbb{N}$, each $a_{j_1,\ldots,j_n} \in \mathbb{C}$, and is zero whenever two or more of the indices $j_1, \ldots, j_n$ coincide and $A_1, \ldots, A_N \in \mathcal{B}(S)$, with $A_i$ of the form $J_i \times B_i$ where $J_i$ is an interval in $[0,T]$ and $B_i \in \mathcal{B}(\mathbb{R})$ with $\rho(B_i) < \infty$, for each $1 \leq i \leq N$. It is shown as in Proposition 1.6 of Huang and Yan [23] that $\mathcal{D}^{(n)}$ is dense in $H^{\otimes n}$. It then follows that $\mathcal{D}_s^{(n)}$ is dense in $H_s^{\otimes n}$ where $\mathcal{D}_s^{(n)} := \mathcal{D}^{(n)} \cap H_s^{\otimes n}$.

For each $f_n \in \mathcal{D}^{(n)}$ we define its *multiple Wiener-Lévy integral* by

$$I_n(f_n) = \sum_{j_1,\ldots,j_n=1}^{N} a_{j_1,\ldots,j_n} M(A_{j_1}) \cdots M(A_{j_n}). \quad (4.25)$$

The mapping $f_n \to I_n(f_n)$ is easily seen to be linear. For each $f_n \in \mathcal{D}^{(n)}$, $I_n(f_n) = I_n(\widehat{f_n})$, where $\widehat{f_n}$ is the symmetrisation of $f_n$.

The next result is due to Itô ([25]). The special case where $M$ is a Brownian motion is proved in many textbooks (see e.g. [39], [23]) and the general case proceeds along similar lines.

**Theorem 4.2** *For each $f_m \in \mathcal{D}_s^{(m)}, g_n \in \mathcal{D}_s^{(n)}, m, n \in \mathbb{N}$*

$$\mathbb{E}(I_m(f_m)) = 0, \quad \mathbb{E}(\overline{I_m(f_m)} I_n(g_n)) = n! \langle f_m, g_n \rangle \delta_{mn}.$$

So for each $n \in \mathbb{N}$, $I_n$ is an isometry from $\mathcal{D}_s^{(n)}$ (equipped with the inner product $\langle\langle \cdot, \cdot \rangle\rangle := n! \langle \cdot, \cdot \rangle$) into $\mathcal{H}$. It hence extends to an isometry which is defined on the whole of $H_s^{\odot n}$. We continue to denote this mapping by $I_n$ and for each $f_n \in H_s^{\odot n}$, we call $I_n(f_n)$ the *multiple Wiener-Lévy integral* of $f_n$. By continuity and Theorem 4.2, we obtain

$$\mathbb{E}(I_m(f_m)) = 0, \quad \mathbb{E}(\overline{I_m(f_m)} I_n(g_n)) = n! \langle f_m, g_n \rangle \delta_{mn}, \quad (4.26)$$

for each $f_m \in H_s^{\odot m}, g_n \in H_s^{\odot n}, m, n \in \mathbb{N}$.

We introduce the $n$-simplex $\Delta_n$ in $[0,T]$ so

$$\Delta_n = \{0 < t_1 < \cdots < t_n < T\}$$

and define the iterated stochastic integral

$$J_n(f_n) := \int_{\Delta_n \times \mathbb{R}^n} f_n(w_1, \ldots, w_n) M(dw_1) \cdots M(dw_n)$$
$$= \int_0^T \int_{\mathbb{R}} \int_0^{t_n-} \int_{\mathbb{R}} \cdots \int_0^{t_2-} \int_{\mathbb{R}} f_n(t_1, x_1, \ldots, t_n, x_n) M(dt_1, dx_1) \cdots M(dt_n, dx_n),$$

for each $f_n \in H^{\otimes n}$. Then if $f_n \in H_s^{\odot n}$, we have

$$I_n(f_n) = n! J_n(f_n). \tag{4.27}$$

This is established by exactly the same argument as the case where $M$ is a Brownian motion (see e.g. [39]) i.e. first establish the result when $f_n \in \mathcal{D}_s^{(n)}$ by a direct (but messy) calculation and then pass to the general case by means of an approximation.

The final result in this section is the celebrated chaos decomposition which is again due to Itô [25].

**Theorem 4.3**
$$\mathcal{H} = \mathbb{C} \oplus \bigoplus_{n=1}^{\infty} Ran(I_n). \tag{4.28}$$

*Proof.* This follows by iteration of the Itô representation as in Theorem 4 of [35]. □

We use $U$ to denote the mapping $\bigoplus_{n=0}^{\infty} I_n$ from $\hat{\Gamma}(H)$ to $\mathcal{H}$, where $I_0(f_0) := f_0$ for all $f_0 \in \mathbb{C}$. $U$ is sometimes called the *duality transform*. The next result dates back to Segal [49] in its Gaussian version. The extension to Lévy processes first seems to have been made explicit by Dermoune [12].

**Corollary 4.1** *[Wiener-Segal-Itô Isomorphism]* $U$ *is a unitary isomorphism between* $\hat{\Gamma}(H)$ *and* $\mathcal{H}$.

*Proof.* It follows from (4.26) that $U$ is an isometry. By (4.28), if $F \in \mathcal{H}$, there exists a sequence $(f_n, n \in \mathbb{N})$ with each $f_n \in H_s^{\odot n}$, such that $F = \sum_{n=0}^{\infty} I_n(f_n)$. Thus we see that $U$ is surjective, and hence is unitary. □

## 4.4 The Role of Stochastic Exponentials

We have seen in section 1 that exponential vectors play an important structural role in Fock space (see also [34], [38], [43]). In this section we will find the analogous vectors in Lévy-Itô space.

Let $Y = (Y(t), 0 \leq t \leq T)$ be a complex valued semimartingale defined on $(\Omega, \mathcal{F}, (\mathcal{F}_t, t \geq 0), P)$ so that each $Y(t) = Y^{(1)}(t) + iY^{(2)}(t)$, where $Y^{(1)}$ and $Y^{(2)}$ are real valued semimartingales. The unique solution to the stochastic differential equation (SDE)

$$dZ(t) = Z(t-)dY(t), \tag{4.29}$$

with initial condition $Z(0) = 1$ (a.s.) is given by the *stochastic exponential* or *Doléans-Dade exponential*,

$$\begin{aligned} Z(t) &= \exp\left\{Y(t) - Y(0) - \frac{1}{2}[Y_c^{(1)}, Y_c^{(1)}](t) + \frac{1}{2}[Y_c^{(2)}, Y_c^{(2)}](t) - i[Y_c^{(1)}, Y_c^{(2)}](t)\right\} \\ &\quad \times \prod_{0 \leq s \leq t} [1 + \Delta Y(s)]e^{-\Delta Y(s)} \end{aligned} \tag{4.30}$$

for each $0 \leq t \leq T$. Details can be found in [17] or [26]. Here $[\cdot, \cdot]$ is the quadratic variation and $Y_c^j$ is the continuous part of $Y^j$ for $j = 1, 2$. Henceforth we will write each $Z(t) := \mathcal{E}_Y(t)$.

For each $f \in H$, we introduce the square-integrable martingales $(Y_f(t), 0 \leq t \leq T)$ defined by

$$Y_f(t) = \int_0^t f(s,x) M(ds, dx) = \sigma \int_0^t f(s) dB(s) + \int_0^t \int_{\mathbb{R}-\{0\}} f(s,x) \tilde{N}(ds, dx),$$

where $f(s) := f(s, 0)$, for all $0 \leq s \leq T$.

**Theorem 4.4** *For all $f \in H$,*

$$Ue(f) = \mathcal{E}_{Y_f}(T).$$

*Proof.* Iterating the SDE (4.29) as in [17], p. 189 and using (4.27) we obtain for each $f \in H$,

$$\begin{aligned} \mathcal{E}_{Y_f}(T) &= 1 + \sum_{n=1}^{\infty} J_n(f^{\otimes n}) \\ &= \sum_{n=0}^{\infty} \frac{1}{n!} I_n(f^{\otimes n}). \end{aligned}$$

Now by Corollary 4.1, we have

$$Uf^{\odot n} = I_n(f^{\otimes n}) \Rightarrow U\left(\frac{f^{\otimes n}}{\sqrt{n!}}\right) = \frac{1}{n!}I_n(f^{\otimes n}),$$

and the result follows on using the strong continuity of $U$. □

**Corollary 4.2** *The linear span of $\{\mathcal{E}_{Y_f}(T), f \in H\}$ is dense in $\mathcal{H}$.*

*Proof.* This follows immediately from Theorem 4.4 and the fact that exponential vectors are total in $\Gamma(H)$. □

**Example 1** Brownian Motion

In this case $X(t) = \sigma B(t)$ for each $t \geq 0$ and it is sufficient to take $S = [0, T]$, so $H = L^2([0, T])$ and each $Y_f(t) = \sigma \int_0^t f(s) dB(s)$. Combining Theorem 4.4 with (4.30) we obtain the well known result:

$$Ue(f) = \exp\left\{\sigma \int_0^T f(s) dB(s) - \frac{\sigma^2}{2} \int_0^T |f(s)^2| ds\right\},$$

see e.g. [43], example 19.9, p.130.

**Example 2** The Poisson Process

Let $(N(t), t \geq 0)$ be a Poisson process having intensity $\lambda > 0$. In this case $\sigma = 0$ and $\nu = \lambda \delta_1$. There is then a unique isomorphism $V$ between $H = L^2(S, \mu)$ and $L^2([0, T])$ given by $(Vf)(s, x) = \frac{f(s,1)}{\sqrt{\lambda}}$, for each $s \in [0, T], x \in \mathbb{R}$. Writing $W := \Gamma(V) U \Gamma(V^{-1})$, we find from (4.30) that for all $g \in L^2([0, T])$, taking $I_g(t) = \frac{1}{\sqrt{\lambda}} \int_0^t g(s) dN(s) - \sqrt{\lambda} \int_0^t g(s) ds = \frac{1}{\sqrt{\lambda}} \sum_{0 \leq s \leq t} g(s) \Delta N(s) - \sqrt{\lambda} \int_0^t g(s) ds$, we obtain

$$We(g) = \exp\left\{-\sqrt{\lambda} \int_0^T g(s) ds\right\} \prod_{0 \leq s \leq T} \left(1 + \frac{\Delta N(s)}{\sqrt{\lambda}} g(s)\right)$$

(c.f. [43], Example 19.11, p.131).

# 5 Malliavin Calculus in Lévy-Itô Space

In this section we will work with operators on the Lévy-Itô space $\mathcal{H}$ which are unitary transforms (in the sense of the duality transform) of those defined in $\Gamma(H)$ in section 1, so $H = L^2(S, \mathcal{B}(S), \mu; \mathbb{C})$ where $S = [0, T] \times \mathbb{R}$. We will frequently use the following result. Let $H_1$ and $H_2$ be complex separable Hilbert spaces and let $V$ and $W$ be unitary isomorphisms between $H_1$ and $H_2$. Let $T$ be a densely defined closed linear operator on $H_1$ with domain $\mathcal{D}$. Then it is easily verified that $VTW^{-1}$ is closed on the dense domain $W\mathcal{D}$.

## 5.1 The Gradient

We define the *gradient* $D$ in $\mathcal{H}$ by $D := (U \otimes I)\nabla^- U^{-1}$ on the domain $U\text{Dom}(\nabla^-)$. Then by Theorems 3.1 and 4.3 we see that

$$\text{Dom}(D) = \left\{ \sum_{n=0}^{\infty} I_n(f_n); \sum_{n=1}^{\infty} nn! ||f_n||^2 < \infty \right\}.$$

The spaces $U\Xi$ and $(U \otimes I)\Xi \otimes H$ are infinite dimensional Sobolev spaces which are usually denoted $\mathbb{D}_1^2$ and $\mathbb{D}_1^2(H)$ (respectively) in this context (see e.g. [23]). For each $f \in H, \psi \in \text{Dom}(\nabla^-), Ua(f)\psi = D_f U\psi$ where $D_f := (I \otimes \epsilon_f)D$ is the directional gradient operator. For each $s = (t, x) \in S$, the Malliavin derivative $D_s$ is obtained by writing $D = (D_s, s \in S)$ and for each $\psi \in \text{Dom}(\nabla^-)$ we have $U\nabla_s^- \psi = D_s U\psi$ ($\mu$ a.e.)

In particular, if $\mathcal{E}_f(T)$ is the stochastic exponential of $f \in H$ then we have

$$D_s \mathcal{E}_f(T) = f(s)\mathcal{E}_f(T),$$

for all $s \in S$ except for a set of $\mu$-measure zero (c.f. [13], formula (I.17)) and from (3.19) we see that if $\psi = \sum_{n=0}^{\infty} I_n(f_n) \in \text{Dom}(D)$, then

$$D_s \psi(\cdot) = \sum_{n=1}^{\infty} n I_{n-1}(f_n(\cdot, s)) \quad \mu \text{ a.e.},$$

where $f_n(\cdot, t)$ is the symmetric function of $n-1$ variables obtained by fixing $s \in S$.

## 5.2 The Divergence

We define the *divergence* $\delta$ in $\mathcal{H} \otimes H = L^2(S, \mathcal{S}, \mu; \mathcal{H})$ by the prescription $\delta = U\nabla^+(U^{-1} \otimes I)$ on the domain $(U \otimes I)\text{Dom}(\nabla^+)$. Elements of $H^{\odot n} \otimes H$ are sequences $(g_{n+1}, n \in \mathbb{Z}_+)$ where each $g_{n+1}$ is a measurable function of $n+1$ variables which is symmetric in the first $n$ of these. We denote the symmetrisation of $g_{n+1}$ by $\widetilde{g_{n+1}}$ so that for each $s_1, \ldots, s_n, s_{n+1} \in S$,

$$\widetilde{g_{n+1}}(s_1, \ldots, s_n, s_{n+1}) = \frac{1}{n+1} \left( g_{n+1}(s_1, \ldots, s_n, s_{n+1}) + \sum_{j=1}^{n} g_{n+1}(s_1, \ldots, s_{j-1}, s_{n+1}, s_{j+1}, \ldots, s_n, s_j) \right).$$

In the following, we will use the standard notation $I_n(g_{n+1}) := (I_n \otimes I)(g_{n+1})$ so the multiple Wiener-Lévy integral only acts on the symmetric part of the function.

We then have the following characterisation of $\text{Dom}(\delta)$.

**Theorem 5.1** $X = \sum_{n=0}^{\infty} I_n(g_{n+1}) \in Dom(\delta)$ if and only if $\sum_{n=0}^{\infty} (n+1)! ||\widetilde{g_{n+1}}||^2 < \infty$. We then have

$$\delta(X) = \sum_{n=0}^{\infty} I_{n+1}(\widetilde{g_{n+1}}).$$

*Proof.* It is sufficient to take $g_{n+1} = f^{\odot n} \otimes h$, then for each $s_1, \ldots, s_n, s_{n+1} \in S$,

$$\text{Symm}(f^{\odot n}, h)(s_1, \ldots, s_n, s_{n+1})$$
$$= \left( \sqrt{\frac{n!}{n+1}} \sum_{r=0}^{n} f^{\otimes n-r} \otimes h \otimes f^{\otimes r} \right)(s_1, \ldots, s_n, s_{n+1})$$
$$= \sqrt{\frac{n!}{n+1}} \left( f(s_1) \cdots f(s_n) h(s_{n+1}) + \sum_{j=1}^{n} f(s_1) \cdots f(s_{j-1}) h(s_j) f(s_{j+1}) \cdots f(s_n) f(s_j) \right)$$
$$= \sqrt{(n+1)!} \widetilde{g_{n+1}}(s_1, \ldots, s_n, s_{n+1}).$$

The result follows easily from here. $\square$

For each $f \in H$ we define the directional divergence operator $\delta_f$ on $U\mathcal{E}$ by $\delta_f := \delta \circ (I \otimes \epsilon_f^{\dagger})$, then $\delta_f \psi = a^{\dagger}(f) U^{-1} \psi$ for all $\psi \in U\mathcal{E}$.

Arguing as in Theorem 3 of [41] (see also Proposition 3.2 of [14]) it follows that $\delta(X)$ is a non-anticipating extension of the Itô integral $\int_S X(w) M(dw)$ which is defined in the case where $X$ is predictable (see e.g. [2]). We may now rewrite equation (3.20) as the well-known Itô-Skorohod isometry:

$$\mathbb{E}(\overline{\delta(X)}\delta(Y)) = \int_S \mathbb{E}(\overline{X(s)}Y(s))\mu(ds) + \int_S \int_S \mathbb{E}(\overline{D_t X(s)} D_s Y(t))\mu(ds)\mu(dt), \tag{5.31}$$

for all $X, Y \in \mathbb{D}_1^2(H)$ (c.f. Theorem 3.14 in [14]).

## 5.3 Independence Structure

For each $f \in H$, write $f = f_1 + f_2$, where $f_1(t, x) := \left\{ \begin{array}{ll} f(t, 0) & \text{if } x = 0 \\ 0 & \text{if } x \neq 0 \end{array} \right\}$ and $f_2(t, x) := \left\{ \begin{array}{ll} 0 & \text{if } x = 0 \\ f(t, x) & \text{if } x \neq 0 \end{array} \right\}$. We thus obtain a canonical isomorphism between $H$ and $H_1 \oplus H_2$ where $H_1 = L^2([0, T], \lambda_\sigma)$ ($\lambda_\sigma := \sigma^2 \lambda$ is rescaled Lebesgue measure) and $H_2 = L^2(E, \lambda \otimes \nu)$ where $E := [0, T] \times (\mathbb{R} - \{0\})$. Now suppose that $(\Omega, \mathcal{F}, P)$ is of the form $(\Omega_1 \times \Omega_2, \mathcal{F}_1 \otimes \mathcal{F}_2, P_1 \times P_2)$. The canonical example of this is called *Lévy space* in [12] and *Wiener-Poisson* *space* in [24]. In this set-up $\Omega_1$ is the space of continuous functions which vanish at zero equipped with Wiener measure $P_1$ on the $\sigma$-algebra $\mathcal{F}_1$ generated by the cylinder sets. $\Omega_2$ is the set of all $\hat{\mathbb{Z}}_+$-valued measures on $E$ (where $\hat{\mathbb{Z}}_+ := \mathbb{Z}_+ \cup \infty$). $\mathcal{F}_2$ is the smallest $\sigma$-algebra of subsets of $\Omega_2$ which permits all evaluations of measures on Borel sets in $E$ to be measurable and $P_2$ is taken to be a Poisson measure on $(\Omega_2, \mathcal{F}_2)$ with intensity $\lambda \times \nu$.

Let $\mathcal{F}_{1,T} := \sigma\{B(s), 0 \leq s \leq T\}$ and $\mathcal{F}_{2,T} := \sigma\left\{\int_{\mathbb{R}-\{0\}} x\tilde{N}(s, dx), 0 \leq s \leq T\right\}$ so $\mathcal{F}_{i,T}$ are sub-$\sigma$-algebras of $\mathcal{F}_i$ for $i = 1, 2$. Applying Corollary 4.1 separately in $H_1$ and $H_2$ we see that the duality transform factorises as $U = U_1 \otimes U_2$ where $U_i$ is the duality transform between $\Gamma(H_i)$ and $L^2(\Omega_i, \mathcal{F}_{i,T}, P_i)$ for $i = 1, 2$. Applying this to the tensor decompositions (3.16) and (3.17) We then obtain
$$D = \pi_1(D_B \otimes I) + \pi_2(I \otimes D_N),$$
on $\text{Dom}(D_B) \underline{\otimes} \text{Dom}(D_N)$, and
$$\delta = (\delta_B \otimes I)\pi_1^* + (I \otimes \delta_N)\pi_2^*,$$
on $\text{Dom}(\delta_B) \oplus \text{Dom}(\delta_N)$. Here $D_B$ and $\delta_B$ are the usual gradient and divergence associated to Brownian motion (see e.g. [23], [36], [39], [50]) while $D_N$ and $\delta_N$ are the gradient and divergence associated to Poisson random measures (see e.g. [14], [35], [40], [41], [45], [46]).

### 5.4 Number Operator

The number operator in Lévy-Itô space is $\mathcal{N} = UNU^{-1}$ and the corresponding semigroup is $(\mathcal{T}_t, t \geq 0)$ where $\mathcal{T}_t = UT_tU^{-1}$ for each $t \geq 0$. We observe that by Theorem 3.2 we have
$$\delta D = \mathcal{N},$$
on $\text{Dom}(\mathcal{N})$ and by Theorem 3.3 we obtain the $L^2$-equivalence of norms whereby the operator $(I + \mathcal{N})^{-\frac{1}{2}}$ is a unitary isomorphism between $\mathcal{H}$ and $\mathbb{D}_1^2$.

If we employ the independence structure we have
$$\mathcal{N} = \mathcal{N}_B \otimes I + I \otimes \mathcal{N}_N,$$
on $\text{Dom}(\mathcal{N}_B) \otimes \text{Dom}(\mathcal{N}_N)$, and
$$\mathcal{T}_t = \mathcal{T}_t^B \otimes \mathcal{T}_t^N,$$
for all $t \geq 0$, where the sub/superscripts $B$ and $N$ refer to the Brownian and Poisson components in the obvious way.

$-\mathcal{N}_B$ is the well-known infinite-dimensional Ornstein-Uhlenbeck operator which enjoys the hypercontractivity property. Surgailis [52] has shown that $-\mathcal{N}_P$ does not have this property. Furthermore in Theorem 5.1 of [53], Surgailis proves that if $(R_t, t \geq 0)$ is a contraction semigroup in $L^2(E, \lambda \times \mu)$ then the contraction semigroup $(\mathcal{R}_t, t \geq 0)$ where each $\mathcal{R}_t = U_2 \Gamma(R_t) U_2^{-1}$ is positivity preserving if and only if $(R_t, t \geq 0)$ is doubly Markovian. In this latter case, $(\mathcal{R}_t, t \geq 0)$ is itself Markovian. So we can assert the Markovianity of $(\mathcal{T}_t^N, t \geq 0)$ and hence of $(\mathcal{T}_t, t \geq 0)$. Further studies of the semigroup $(\mathcal{T}_t^N, t \geq 0)$ can be found in [10] and [54].

We finish this article with a word of warning. The *universal Malliavin calculus* that we have described here shows great promise for obtaining more wideranging properties which hold generally for processes which enjoy a chaos decomposition. However there will be local features of the calculus which are particular to the process under consideration and which cannot be obtained through the Fock space isomorphism. For an example, see formula (I.25) in [13] for an algebraic relation between the divergence and the gradient where there is an extra term in the Poisson case which is absent in Gaussian spaces.

**Acknowledgements.** This paper grew out of two lecture courses which I gave at the Universities of Sheffield (2005) and Virginia (2006), respectively. I would like to thank all the participants in these for their contributions. Particular thanks are due to Len Scott for arranging the trip to Virginia and for his warm and generous hospitality during my visit. I am also grateful to Martin Lindsay for a helpful discussion about $\nabla^-$ and $\nabla^+$ and to Nick Bingham for drawing my attention to [9]. I would like to thank Fangjun Hsu for reading through an early draft of this article and making some very helpful suggestions. Thanks are also due to both Nick Bingham and Robin Hudson for valuable comments in this regard.